\documentclass{amsart}

\usepackage{amsmath,amsthm,amsfonts,amscd,amssymb}

\newtheorem{thm}{Theorem}[section]
\newtheorem{cor}[thm]{Corollary}
\newtheorem{lem}[thm]{Lemma}
\newtheorem{prop}[thm]{Proposition}

\newtheorem{prob}{Problem}

\theoremstyle{definition}

\theoremstyle{remark}

\begin{document}

\title[Search bounds for zeros of polynomials over $\overline{\mathbb Q}$]{Search bounds for zeros of polynomials over the algebraic closure of $\mathbb Q$}
\author{Lenny Fukshansky}

\address{Department of Mathematics, Mailstop 3368, Texas A\&M University, College Station, Texas 77843-3368}
\email{lenny@math.tamu.edu}

\begin{abstract}
We discuss existence of explicit search bounds for zeros of polynomials with coefficients in a number field. Our main result is a theorem about the existence of polynomial zeros of small height over the field of algebraic numbers outside of unions of subspaces. All bounds on the height are explicit. 
\end{abstract}

\maketitle

\def\A{{\mathcal A}}
\def\B{{\mathcal B}}
\def\C{{\mathcal C}}
\def\D{{\mathcal D}}
\def\F{{\mathcal F}}
\def\x{{\mathcal H}}
\def\I{{\mathcal I}}
\def\J{{\mathcal J}}
\def\K{{\mathcal K}}
\def\L{{\mathcal L}}
\def\M{{\mathcal M}}
\def\R{{\mathcal R}}
\def\s{{\mathcal S}}
\def\V{{\mathcal V}}
\def\X{{\mathcal X}}
\def\Y{{\mathcal Y}}
\def\H{{\mathcal H}}
\def\cee{{\mathbb C}}
\def\Nn{{\mathbb N}}
\def\pee{{\mathbb P}}
\def\que{{\mathbb Q}}
\def\real{{\mathbb R}}
\def\zed{{\mathbb Z}}
\def\gmn{{\mathbb G_m^N}}
\def\qbar{{\overline{\mathbb Q}}}
\def\eps{{\varepsilon}}
\def\vek{{\varepsilon_k}}
\def\ahat{{\hat \alpha}}
\def\bhat{{\hat \beta}}
\def\gt{{\tilde \gamma}}
\def\h{{\tfrac12}}
\def\ba{{\boldsymbol a}}
\def\be{{\boldsymbol e}}
\def\bei{{\boldsymbol e_i}}
\def\bc{{\boldsymbol c}}
\def\bm{{\boldsymbol m}}
\def\bk{{\boldsymbol k}}
\def\bi{{\boldsymbol i}}
\def\bl{{\boldsymbol l}}
\def\bq{{\boldsymbol q}}
\def\bu{{\boldsymbol u}}
\def\bt{{\boldsymbol t}}
\def\bs{{\boldsymbol s}}
\def\bv{{\boldsymbol v}}
\def\bw{{\boldsymbol w}}
\def\bx{{\boldsymbol x}}
\def\bX{{\boldsymbol X}}
\def\bz{{\boldsymbol z}}
\def\bwy{{\boldsymbol y}}
\def\bg{{\boldsymbol g}}
\def\bY{{\boldsymbol Y}}
\def\bL{{\boldsymbol L}}
\def\baa{{\boldsymbol\alpha}}
\def\bb{{\boldsymbol\beta}}
\def\bet{{\boldsymbol\eta}}
\def\bxi{{\boldsymbol\xi}}
\def\bo{{\boldsymbol 0}}
\def\bol{{\boldsymbol 1}_L}
\def\ep{\varepsilon}
\def\p{\boldsymbol\varphi}
\def\q{\boldsymbol\psi}
\def\rank{\operatorname{rank}}
\def\aut{\operatorname{Aut}}
\def\lcm{\operatorname{lcm}}
\def\sgn{\operatorname{sgn}}
\def\spn{\operatorname{span}}
\def\md{\operatorname{mod}}
\def\Norm{\operatorname{Norm}}
\def\dim{\operatorname{dim}}
\def\det{\operatorname{det}}
\def\Vol{\operatorname{Vol}}
\def\rk{\operatorname{rk}}

\section{Introduction}
Let $F_1,...,F_k$ be a collection of nonzero polynomials in $N$ variables of respective degrees $M_1,...,M_k$ with coefficients in a number field $K$ of degree $d$ over $\que$. Consider a system of equations
\begin{equation}
\label{i1}
F_1(X_1,...,X_N) = ... = F_k(X_1,...,X_N) = 0.
\end{equation}
There are two fundamental questions one can ask about this system: does (\ref{i1}) have nonzero solutions over $K$, and, if yes, how do we find them? In \cite{masser:baker}, D. W. Masser poses these general questions for a system of equations with integer coefficients, and suggests an alternative approach to both of them simultaneously by introducing {\it search bounds} for solutions. We start by generalizing this approach over $K$. 

We write $\qbar$ for the algebraic closure of $\que$, and write $\pee(\qbar^N)$ for the projective space over $\qbar^N$. If $H$ is a height function defined over $\qbar$, then by Northcott's theorem \cite{northcott:ht} a set of the form 
\begin{equation}
\label{i2}
S_D(C) = \{ \bx \in \pee(\qbar^N) : H(\bx) \leq C,\  \deg(\bx) \leq D \}
\end{equation}
has finite cardinality for any $C,D \in \real$, where $\deg(\bx)$ is degree of the field extension generated by the coordinates of $\bx$ over $\que$. Suppose that we were able to prove that if (\ref{i1}) has a nonzero solution $\bx \in K^N$, then it has such a solution with $H(\bx) \leq C$ for some explicit $C$. This means that we can restrict the search for a solution to a subset of the finite set $S_d(C)$ as in (\ref{i2}). We will call a constant $C$ like this a {\it search bound} for (\ref{i1}). If a search bound like this exists, it will clearly depend on heights of polynomials $F_1,...,F_k$. As in \cite{masser:baker}, we can now replace the two questions above by the following problem.

\begin{prob} \label{fund_prob} Find an explicit search bound for a nonzero solution of (\ref{i1}) over~$K$.
\end{prob}

This problem has been solved for arbitrary $N$ only in very few cases. First suppose that $k<N$, and $M_1 = ... = M_k = 1$. If $F_1,...,F_k$ are homogeneous, a solution to Problem \ref{fund_prob} is provided by Siegel's Lemma (see \cite{vaaler:siegel}). In the case when $F_1,...,F_k$ are inhomogeneous linear polynomials, this problem has been solved in \cite{vaaler:oleary}. Another instance of (\ref{i1}) for which the general solution to Problem \ref{fund_prob} is known is that of one quadratic polynomial. If $k=1$, $M_1=2$, and $F_1$ is a quadratic form in $N \geq 2$ variables with coefficients in $K$, a solution to Problem \ref{fund_prob} is presented in \cite{cassels:small} in case $K=\que$, and generalized to an arbitrary number field in \cite{raghavan}. If $F_1$ is an inhomogeneous quadratic polynomial, a general solution to Problem \ref{fund_prob} over $\que$ can be found in \cite{masser:1}, and its generalization to an arbitrary number field in \cite{me:smallzeros}. For a review of further advances in this subject and a detailed bibliography, see \cite{masser:baker}.

A general solution to Problem \ref{fund_prob} even for one polynomial of arbitrary degree in an arbitrary number of variables seems to be completely out of reach at the present time. In fact, if $K=\que$ and $F_1,...,F_k$ are homogeneous, a solution to Problem 1 would provide an algorithm to decide whether a system of homogeneous Diophantine equations has an integral solution, and so would imply a positive answer to Hilbert's 10th problem in this case. However, by the famous theorem of Matijasevich \cite{mat} Hilbert's 10th problem is undecidable. This means that in general search bounds do not exist over $\que$; in fact, they are unlikely to exist over any fixed number field. Moreover, it is known they do not exist over $\que$ even for a single quartic polynomial or for a system of quadratics (see \cite{masser:baker} for details).

In this paper we deal with the case of a single polynomial. Let us relax the condition that a solution must lie over a fixed number field $K$, but instead search for a solution of bounded height and bounded degree over $\qbar$. In other words, given an equation of the form
$$F(X_1,...,X_N) = 0,$$
we want to prove the existence of a nonzero solution $\bx \in \qbar^N$ such that $H(\bx) \leq C$ and $\deg_K(\bx) \leq D$ for explicit constants $C$ and $D$, where  $\deg_K(\bx)$ stands for the degree of the field extension over $K$ generated by the coordinates of $\bx$. This problem is easily tractable as we will show in section~3, and still provides an explicit search bound since the set $S_D(C)$ is finite. In fact, we can prove a stronger statement by requiring the point $\bx$ in question to satisfy some additional arithmetic conditions. Write $\gmn$ for the multiplicative torus $(\qbar^{\times} )^N$. Here is the main result of this paper.

\begin{thm} \label{main} Let $F(X_1,...,X_N)$ be a homogeneous polynomial in $N \geq 2$ variables of degree $M \geq 1$ over a number field $K$, and let $A \in GL_N(K)$. Then either there exists $\boldsymbol 0 \neq \bx \in K^N$ such that $F(\bx)=0$ and
\begin{equation}
\label{main:1}
H(\bx) \leq H(A),
\end{equation}
or there exists $\bx \in A\gmn$ with $\deg_K(\bx) \leq M$ such that $F(\bx) = 0$, and
\begin{equation}
\label{main:2}
H(\bx) \leq C_1(N,M) H(A)^2 H(F)^{1/M},
\end{equation}
where
\begin{equation}
\label{main:3}
C_1(N,M) = 2^{N-1} \left( \frac{M+2}{2} \right)^{\frac{(4M+1)(N-2)}{2M}} \binom{M+N}{N}^{\frac{1}{2M}} \prod_{j=2}^N \binom{M+j-2}{j-2}^{\frac{1}{2M}}.
\end{equation}
\end{thm}

\noindent
In other words, Theorem \ref{main} asserts that for each element $A$ of $GL_N(K)$ either there exists a zero of $F$ over $K$ whose height is bounded by $H(A)$, or there exists a small-height zero of $F$ over $\qbar$ which lies outside of the union of nullspaces of row vectors of $A^{-1}$; for instance, if $A=I_N$ this means that there exists a small-height zero of $F$ with all coordinates non-zero.

 Notice that our approach of searching for small-height polynomial zeros over $\qbar$ is analogous in spirit to the so called ``absolute'' results, like the absolute Siegel's Lemma of Roy and Thunder, \cite{absolute:siegel}. The difference however is that we also keep a bound on the degree of a solution over the base field $K$. 

This paper is organized as follows. In section~2 we set the notation and introduce the height functions that we will use. In section~3 we talk about the basic search bounds for zeros of a given polynomial over $\qbar$. In section~4 we prove Theorem \ref{main}. Results of this paper also appear as a part of \cite{me:diss}.

\section{Notation and heights}
We start with some notation. Let $K$ be a number field of degree $d$ over $\que$, $O_K$ its ring of integers, and $M(K)$ its set of places. For each place $v \in M(K)$ we write $K_v$ for the completion of $K$ at $v$ and let $d_v = [K_v:\que_v]$ be the local degree of $K$ at $v$, so that for each $u \in M(\que)$
\begin{equation}
\sum_{v \in M(K), v|u} d_v = d.
\end{equation}

\noindent
For each place $v \in M(K)$ we define the absolute value $\|\ \|_v$ to be the unique absolute value on $K_v$ that extends either the usual absolute value on $\real$ or $\cee$ if $v | \infty$, or the usual $p$-adic absolute value on $\que_p$ if $v|p$, where $p$ is a prime. We also define the second absolute value $|\ |_v$ for each place $v$ by $|a|_v = \|a\|_v^{d_v/d}$ for all $a \in K$. Then for each non-zero $a \in K$ the {\it product formula} reads
\begin{equation}
\label{product_formula}
\prod_{v \in M(K)} |a|_v = 1.
\end{equation} 

\noindent
For each $v \in M(K)$ define a local height $H_v$ on $K_v^N$ by
\[ H_v(\bx) = \left\{ \begin{array}{ll}
\max_{1 \leq i \leq N} |x_i|_v & \mbox{if $v \nmid \infty$} \\
\left( \sum_{i=1}^N \|x_i\|_v^2 \right)^{d_v/2d} & \mbox{if $v | \infty$}
\end{array}
\right. \]
for each $\bx \in K_v^N$. We define the following global height function on $K^N$:
\begin{equation}
H(\bx) = \prod_{v \in M(K)} H_v(\bx),
\end{equation}
for each $\bx \in K^N$. Notice that due to the normalizing exponent $1/d$, our global height function is absolute, i.e. for points over $\qbar$ its value does not depend on the field of definition. This means that if $\bx \in \qbar^N$ then $H(\bx)$ can be evaluated over any number field containing the coordinates of $\bx$.

We also define a height function on algebraic numbers. Let $\alpha \in \qbar$, and let $K$ be a number field containing $\alpha$. Then define 
 \begin{equation}
h(\alpha) = \prod_{v \in M(K)} \max\{1, |\alpha|_v\}.
\end{equation}

\noindent
We define the height of a polynomial to be the height of the corresponding coefficient vector. We also define height on $GL_N(K)$ by viewing matrices as vectors in $K^{N^2}$. On the other hand, if $M<N$ are positive integers and $A$ is an $M \times N$ matrix with row vectors $\ba_1,...,\ba_M$, we let
\begin{equation}
H(A) = H(\ba_1 \wedge ... \wedge \ba_M),
\end{equation}
and if $V$ is the nullspace of $A$ over $K$, we define $H(V) = H(A)$. This is well defined, since multiplication by an element of $GL_M(K)$ does not change the height. In other words, for a subspace $V$ of $K^N$ its height is defined to be the height of the corresponding point on a Grassmannian. 

We will need the following basic property of heights, which can be easily derived from Lemma 2 of \cite{vaaler:pinner} (see Lemma 4.1.1 of \cite{me:diss} for details).

\begin{lem} \label{mahler} Let $g(X) \in K[X]$ be a polynomial of degree $M$ in one variable with coefficients in $K$. There exists $\alpha \in \qbar$ of degree at most $M$ over $K$ such that $g(\alpha) = 0$, and
\begin{equation}
\label{bz3}
h(\alpha) \leq H(g)^{1/M}.
\end{equation}
\end{lem}
\smallskip

Throughout this paper, let $M,N$ be positive integers, and define
\begin{equation}
\label{bz4}
\M(N,M) = \left\{ (i_1,...,i_N) \in \zed_+^N : \sum_{j=1}^N i_j = M \right\},
\end{equation}
where $\zed_+$ is the set of all non-negative integers. Then any homogeneous polynomial $F$ in $N$ variables of degree $M$ with coefficients in $K$ can be written as
$$F(X_1,...,X_N) = \sum_{\bi \in \M(N,M)} f_{\bi} X_1^{i_1} \dots X_N^{i_N} \in K[X_1,...,X_N].$$

For a point $\bz = (z_1,...,z_N) \in \qbar^N$, we write $\deg_K(\bz)$ to mean the degree of the extension $K(z_1,...,z_N)$ over $K$, i.e. $\deg_K(\bz) = [K(z_1,...,z_N):K]$. We are now ready to proceed.

\section{Basic bounds for one polynomial}

We start by exhibiting a basic bound for zeros of polynomials over $\qbar$.

\begin{prop} \label{bz:basic} Let $M \geq 1$, $N \geq 2$, and $F(X_1,...,X_N)$ be a homogeneous polynomial in $N$ variables of degree $M$ with coefficients in a number field $K$. There exists $\boldsymbol 0 \neq \bz \in \qbar^N$ with $\deg_K(\bz) \leq M$ such that $F(\bz) = 0$ and
\begin{equation}
\label{bz:basic_bound}
H(\bz) \leq \sqrt{2}\ H(F)^{1/M}.
\end{equation}
\end{prop}

\proof
If $F$ is identically zero, then we are done. So assume $F$ is non-zero. Write $\be_1,...,\be_N$ for the standard basis vectors for $\qbar^N$ over $\qbar$. Assume that for some $1 \leq i \leq N$, $\deg_{X_i} F < M$, then it is easy to see that $F(\be_i)=0$, and $H(\be_i) = 1$. If $N>2$, let
$$F_1(X_1,X_2) = F(X_1,X_2,0,...,0),$$
and a point $\bx = (x_1,x_2) \in \qbar^2$ is a zero of $F_1$ if and only if $(x_1,x_2,0,...,0)$ is a zero of $F$, and $H(x_1,x_2) = H(x_1,x_2,0,...,0)$. In particular, if $F_1(X_1,X_2) = 0$, then $F(\be_1) = 0$. Hence we can assume that $N=2$, $F(X_1,X_2) \neq 0$, and $\deg_{X_1} F = \deg_{X_2} F = M$. Write
$$F(X_1,X_2) = \sum_{i=0}^M f_i X_1^i X_2^{M-i},$$
where $f_0, f_M \neq 0$. Let
$$g(X_1) = F(X_1,1) = \sum_{i=0}^M f_i X_1^i \in K[X_1],$$
be a polynomial in one variable of degree $M$ with coefficients in $K$. Notice that since coefficients of $g$ are those of $F$, we have $H(g) = H(F)$. By Lemma \ref{mahler}, there must exist $\alpha \in \qbar$ with $\deg_K(\alpha) \leq M$ such that $g(\alpha) = 0$, and
$$H(\alpha,1) \leq \sqrt{2}\ h(\alpha) \leq \sqrt{2}\ H(g)^{1/M} = \sqrt{2}\ H(F)^{1/M}.$$
Taking $\bz = (\alpha,1)$, completes the proof.
\endproof

Notice that if $N=2$, then the bound (\ref{bz:basic_bound}) is best possible with respect to the exponent. Take
$$F(X_1,X_2) = X_1^M - C X_2^M,$$
for some $0 \neq C \in K$. Then zeros of $F$ are of the form $(\alpha C^{1/M}, \alpha)$ for $\alpha \in \qbar$, and it is easy to see that $H(\alpha C^{1/M}, \alpha) \geq \frac{1}{\sqrt{2}} H(F)^{1/M}$.

\begin{cor} Let the notation be as in Proposition \ref{bz:basic}. Then there exist vectors $\bx_{ij} \in \qbar^N$ with non-zero coordinates $i$-th and $j$-th coordinates, $1 \leq i \neq j \leq N$, and the rest of the coordinates equal to zero such that $F(\bx_{ij}) = 0$, $\deg_K(\bx_{ij}) \leq M$, and each $\bx_{ij}$ satisfies (\ref{bz:basic_bound}). Notice that $\qbar^N = \spn_{\qbar} \{ \bx_{ij} : 1 \leq i \neq j \leq N \}$. 
\end{cor}

\proof
In the proof of Proposition \ref{bz:basic} instead of setting all but $X_1$ and $X_2$ equal to zero, set all but $X_i$ and $X_j$ equal to zero. 
\endproof

\section{Proof of Theorem \ref{main}}

Notice that Proposition \ref{bz:basic} only proves the existence of a small-height zero of $F$ which is {\it degenerate} in the sense that it really is a zero of a binary form to which $F$ is trivially reduced. Do there necessarily exist {\it non-degenerate} zeros of $F$? To answer this question, we consider the problem of Proposition \ref{bz:basic} with additional arithmetic conditions. We wonder what can be said about zeros of a polynomial over $\qbar$ outside of a collection of subspaces? For instance, under which conditions does a polynomial $F$ vanish at a point with nonzero coordinates? Here is a simple effective criterion.

\begin{prop} \label{criterion} Let $N \geq 2$, and let $F(X_1,...,X_N) \in K[X_1,...,X_N]$ have degree $M \geq 1$. If $F$ is not a monomial, then there exists $\bz \in \qbar^N$ with $\deg_K(\bz) \leq M$ such that $F(\bz)=0$, $z_i \neq 0$ for all $1 \leq i \leq N$, and
\begin{equation}
\label{nz}
H(\bz) \leq M^M \sqrt{N-1} H(F).
\end{equation}
\end{prop}

\proof
Since $F$ is not a monomial, there must exist a variable which is present to different powers in at least two different monomials, we can assume without loss of generality that it is $X_1$. Then we can write
$$F(X_1,...,X_N) = \sum_{i=0}^M F_i(X_2,...,X_N) X_1^i,$$
where each $F_i$ is a polynomial in $N-1$ variables of degree at most $M-i$. At least two of these polynomials are not identically zero, say $F_j$ and $F_k$ for some $0 \leq j < k \leq M$. Let
$$F_{jk}(X_2,...,X_N) = F_j(X_2,...,X_N) F_k(X_2,...,X_N),$$
then $F_{jk}$ has degree at most $2M-1$. By Lemma 2.2 of \cite{me:classical}, there exists $\ba \in \zed^{N-1}$ such that $a_i \neq 0$ for all $2 \leq i \leq N-1$, $F_{jk}(\ba) \neq 0$, and 
$$\max_{1 \leq i \leq N-1} |a_i| \leq M,$$
hence $H(\ba) \leq M \sqrt{N-1}$. Then $g(X_1) = F(X_1,a_2,...,a_N)$ is a polynomial in one variable of degree at most $M$ with at least two nonzero monomials. If $v \in M(K)$, $v \nmid \infty$, then $H_v(g) \leq H_v(F)$. If $v | \infty$, then for each $0 \leq i \leq M$ we have $\| F_i(\ba) \|_v \leq M^{M-i} H_v(F_i)$, and so
\begin{equation}
\label{g_bound}
H(g) \leq M^{M-1} H(F).
\end{equation}
By factoring a power of $X_1$, if necessary, we can assume that $g$ is a polynomial of degree at least one with coefficients in $K$ such that $g(0) \neq 0$. Then, combining Lemma \ref{mahler} with (\ref{g_bound}), we see that there exists $0 \neq \alpha \in \qbar$ such that $[K(\alpha):K] \leq M$, $g(\alpha) = 0$, and
$$h(\alpha) \leq H(g) \leq M^{M-1} H(F).$$
Let $\bz = (\alpha,\ba)$, then $F(\bz) = 0$, $\deg_K(\bz) \leq M$, $z_i \neq 0$ for each $1 \leq i \leq N$, and
$$H(\bz) \leq h(\alpha) H(\ba) \leq M^M \sqrt{N-1} H(F).$$
\endproof

Under stronger conditions we can find a zero of $F$ of smaller height, all coordinates of which are non-zero.

\begin{thm} \label{nonzero:coord} Let $F(X_1,...,X_N)$ be a homogeneous polynomial in $N \geq 2$ variables of degree $M \geq 1$ with coefficients in a number field $K$. Suppose that $F$ does not vanish at any of the standard basis vectors $\be_1,...,\be_N$. Then there exists $\bz \in \qbar^N$ with $\deg_K(\bz) \leq M$ such that $F(\bz) = 0$, $z_i \neq 0$ for all $1 \leq i \leq N$, and
\begin{equation}
\label{coord:bound}
H(\bz) \leq C_2(N,M)\ H(F)^{1/M},
\end{equation}
where 
\begin{equation}
\label{coord:constant}
C_2(N,M) = 2^{N-1} \left( \frac{M+2}{2} \right)^{\frac{(4M+1)(N-2)}{2M}} \prod_{j=2}^N \binom{M+j-2}{j-2}^{\frac{1}{2M}}.
\end{equation}
\end{thm}

\proof
We argue by induction on $N$. If $N=2$, then the result follows from the argument in the proof of Proposition \ref{bz:basic}. Assume $N>2$. Let $\beta$ be a positive integer, and let
$$F'_{\pm \beta}(X_1,...,X_{N-1}) = F(X_1,...,X_{N-1},\pm \beta X_{N-1}),$$
in other words set $X_N=\pm \beta X_{N-1}$, where the choice of $\pm \beta$ is to be specified later. Let $\be'_1,...,\be'_{N-1}$ be the standard basis vectors for $\qbar^{N-1}$. Notice that if $F'_{\pm \beta}$ vanishes at $\be'_i$ for $1 \leq i \leq N-2$, then $F$ vanishes at $\be_i$, which is a contradiction. In particular, $F'_{\pm \beta}$ cannot be a monomial and cannot be identically zero. Suppose that $F'_{\pm \beta}(\be'_{N-1})=0$. This means that $F'_{\pm \beta}(0,...,0,X_{N-1})$ is identically zero. Write $\bu_i = (0,...,0,i,M-i) \in \zed^N$ for each $0 \leq i \leq M$. Let
$$G(X_{N-1},X_N) = F(0,...,0,X_{N-1},X_N) = \sum_{i=0}^M f_{\bu_i} X_{N-1}^i X_N^{M-i},$$
then
$$F'_{\pm \beta}(0,...,0,X_{N-1}) = G(X_{N-1},\pm \beta X_{N-1}) = \left( \sum_{i=0}^M f_{\bu_i} (\pm \beta)^{M-i} \right) X_{N-1}^M = 0,$$
that is
\begin{equation}
\label{bz7}
\sum_{i=0}^{M} f_{\bu_i} (\pm \beta)^{M-i} = 0.
\end{equation}
Notice that $f_{\bu_0} \neq 0$ and $f_{\bu_M} \neq 0$,  since otherwise $F(\be_N) = 0$ or $F(\be_{N-1}) = 0$. Therefore the left hand side of (\ref{bz7}) is a non-zero polynomial of degree $M$ in $\beta$, and $0$ is not one of its roots, so it has $M$ non-zero roots. Therefore for the appropriate choice of $\pm$ we can select $\beta \in \zed_+$ such that (\ref{bz7}) is {\it not} true and 
\begin{equation}
\label{bz8}
0 < \beta \leq \frac{M}{2}+1 = \frac{M+2}{2}.
\end{equation}
Then for this choice of $\pm \beta$, $F'_{\pm \beta}$ is a polynomial in $N-1$ variables of degree $M$ which does not vanish at any of the standard basis vectors. From now on we will write $F'_{\beta}$ instead of $F'_{\pm \beta}$ for this fixed choice of $\pm \beta$.

Next we want to estimate height of such $F'_{\beta}$. Let $\bl \in \zed_+^{N-1}$ be such that $\sum_{i=1}^{N-1} l_i = M$. There exist $l_{N-1}+1 \leq M+1$ vectors $\bm_j \in \zed_+^N$ such that $m_{ji} = l_i$ for each $1 \leq i \leq N-2$ and $m_{j(N-1)}+m_{jN} = l_{N-1}$, where $0 \leq j \leq l_{N-1}$. Therefore the monomial of $F'_{\beta}$ which is indexed by $\bl$ will have coefficient 
\begin{equation}
\label{bz9}
\alpha_{\bl} = \sum_{j=0}^{l_{N-1}} f_{\bm_j} (\pm \beta)^{l_{N-1}-j}.
\end{equation}
Then for each $v \nmid \infty$
\begin{equation}
\label{bz10}
|\alpha_{\bl}|_v \leq H_v(F),
\end{equation}
and for each $v | \infty$
\begin{eqnarray}
\label{bz11}
\|\alpha_{\bl}\|_v^2 & \leq & \sum_{i=0}^{l_{N-1}} \sum_{j=0}^{l_{N-1}} \beta^{2l_{N-1}-i-j} \|f_{\bm_i}\|_v \|f_{\bm_j}\|_v \nonumber \\
& \leq & \left( \frac{\beta^{2l_{N-1}}}{2} \right) \sum_{i=0}^{l_{N-1}} \sum_{j=0}^{l_{N-1}} (\|f_{\bm_i}\|^2_v + \|f_{\bm_j}\|^2_v) \nonumber \\
& \leq & \left( \frac{\beta^{2l_{N-1}} (l_{N-1}+1)}{2} \right) \left( \sum_{i=0}^{l_{N-1}} \|f_{\bm_i}\|^2_v + \sum_{j=0}^{l_{N-1}} \|f_{\bm_j}\|^2_v \right) \nonumber \\
& \leq & \beta^{2M} (M+2) H_v(F)^2 \leq 2 \left( \frac{M+2}{2} \right)^{2M+1} H_v(F)^2,
\end{eqnarray}
where the last inequality follows by (\ref{bz8}). Therefore, by (\ref{bz10}) and (\ref{bz11}), we have for each $v \nmid \infty$,
\begin{equation}
\label{bz12}
H_v(F'_{\beta}) \leq H_v(F),
\end{equation}
and for each $v | \infty$,
\begin{eqnarray}
\label{bz13}
H_v(F'_{\beta}) & = & \left( \sum_{\bl \in \M(N-1,M)} \|\alpha_{\bl}\|_v^2 \right)^{1/2} \nonumber \\
& \leq & \sqrt{2}\ |\M(N-1,M)|^{1/2} \left( \frac{M+2}{2} \right)^{\frac{2M+1}{2}} H_v(F) \nonumber \\
& \leq & \sqrt{2}\ \binom{M+N-2}{N-2}^{1/2} \left( \frac{M+2}{2} \right)^{\frac{2M+1}{2}} H_v(F)
\end{eqnarray}
Putting (\ref{bz12}) and (\ref{bz13}) together implies that 
\begin{equation}
\label{bz14}
H(F'_{\beta}) \leq \sqrt{2}\ \binom{M+N-2}{N-2}^{1/2} \left( \frac{M+2}{2} \right)^{\frac{2M+1}{2}} H(F).
\end{equation}
By induction hypothesis, there exists $\bx \in \qbar^{N-1}$ with $\deg_K(\bx) \leq M$ such that $F'_{\beta}(\bx)=0$, $x_i \neq 0$ for all $1 \leq i \leq N-1$,  and
\begin{eqnarray}
\label{bz15}
H(\bx) & \leq & C_2(N-1,M)\ H(F'_{\beta})^{\frac{1}{M}} \nonumber \\
& \leq & C_2(N-1,M)\ 2^{\frac{1}{2M}}\ \binom{M+N-2}{N-2}^{\frac{1}{2M}} \left( \frac{M+2}{2} \right)^{\frac{2M+1}{2M}} H(F)^{\frac{1}{M}}.
\end{eqnarray}
Let $E=K(x_1,...,x_{N-1})$. Set $\bz = (\bx, \pm \beta x_{N-1}) \in E^N$, then $\deg_K(\bz) = [E:K] \leq M$, $F(\bz)=0$, $z_i \neq 0$ for all $1 \leq i \leq N$, and applying (\ref{bz8}) and (\ref{bz15}) we have
\begin{eqnarray}
\label{bz16}
H(\bz) & \leq & \prod_{v \nmid \infty} H_v(\bx) \times \prod_{v | \infty} \left( \beta^2 \|x_{N-1}\|_v^2 + H_v(\bx)^2 \right)^{\frac{d'_v}{2d'}} \leq \sqrt{\beta^2 + 1}\ H(\bx) \nonumber \\
& \leq & 2^{\frac{M+1}{2M}}\ \binom{M+N-2}{N-2}^{\frac{1}{2M}} \left( \frac{M+2}{2} \right)^{\frac{4M+1}{2M}} C_2(N-1,M)\ H(F)^{\frac{1}{M}},
\end{eqnarray}
where the product in (\ref{bz16}) is taken over all places in $M(E)$, and $d'_v$, $d'$ stand for local and global degrees of $E$ over $\que$ respectively. The result follows.
\endproof

\noindent
{\it Proof of Theorem \ref{main}.}
Let $K[\bX]_M$ be the space of homogeneous polynomials of degree $M$ in $N$ variables over $K$. For an element $A \in GL_N(K)$ define a map $\rho_A\ : K[\bX]_M \longrightarrow K[\bX]_M$ (compare with \cite{poorten:vaaler}), given by $\rho_A(F)(\bX) = F(A \bX)$ for each $F \in K[\bX]_M$. It is easy to see that the map $A \longmapsto \rho_A$ is a representation of $GL_N(K)$ in $GL(K[\bX]_M)$.

With notation as in the statement of the theorem, let $G(\bX) = \rho_A(F)(\bX)$. First suppose that $G(\be_i) = F(A \be_i) = 0$ for some $1 \leq i \leq N$. Since $\boldsymbol 0 \neq \bwy = A \be_i \in K^N$ is a row of $A$, it is easy to see that 
$$H(\bwy) \leq H(A),$$
which is (\ref{main:1}). Next assume that $G(\be_i) \neq 0$ for each $1 \leq i \leq N$. By Theorem \ref{nonzero:coord}, there exists $\bz \in \gmn$ such that $G(\bz)=0$, $\deg_K(\bz) \leq M$, and 
$$H(\bz) \leq  C_2(N,M)\ H(G)^{1/M}.$$
Then $\bx = A \bz$ is such that $F(\bx) = 0$, $\deg_K(\bx) \leq M$, and $\bx = A \bz \in A \gmn$. It is easy to see that
\begin{equation}
\label{bz19}
H(\bx) \leq H(A) H(\bz) \leq C_2(N,M)\ H(A) H(G)^{1/M}.
\end{equation} 
We now want to estimate $H(G)$. Let $v \in M(K)$. If $v \nmid \infty$, then
\begin{equation}
\label{blue1}
H_v(G) \leq H_v(A)^M H_v(F),
\end{equation}
and if $v | \infty$, then
\begin{equation}
\label{blue2}
H_v(G) \leq \binom{N+M}{N}^{d_v/2d} H_v(A)^M H_v(F).
\end{equation}
These bounds on local heights are well-known. Essentially identical estimates for a bihomogeneous polynomial in two pairs of variables follow from Lemmas 6, 7, and formula (2.2) of \cite{poorten:vaaler}. The proofs of (\ref{blue1}) and (\ref{blue2}) are similar to the proofs of Lemmas 6 and 7 of \cite{poorten:vaaler}, so we do not include them here to maintain the brevity of exposition. Combining (\ref{blue1}) and (\ref{blue2}), we obtain
\begin{equation}
\label{bz20}
H(G) \leq \binom{N+M}{N}^{1/2} H(A)^M H(F).
\end{equation}
The result follows by combining (\ref{bz19}) and (\ref{bz20}).
\boxed{ }

\begin{cor} \label{inhom} Let $F(X_1,...,X_N) \in K[X_1,...,X_N]$ be an inhomogeneous polynomial of degree $M \geq 1$, $N \geq 2$. Suppose that $F$ does not vanish at any of the standard basis vectors $\be_1,...,\be_N$. Then there exists $\bz \in \qbar^N$ with $\deg_K(\bz) \leq M$ such that $F(\bz) = 0$, $z_i \neq 0$ for all $1 \leq i \leq N$, and 
\begin{equation}
\label{inh:bound}
H(\bz) \leq C_2(N+1,M)\ H(F)^{1/M},
\end{equation}
where the constant $C_2(N+1,M)$ is defined by (\ref{coord:constant}) of Theorem \ref{nonzero:coord}.
\end{cor}

\proof
Homogenize $F$ using the variable $X_0$ and denote the resulting homogeneous polynomial in $N+1$ variables by $F'(X_0,...,X_N)$. Then $F'$ has degree $M$, its coefficients are in $K$, and
$$F(X_1,...,X_N) = F'(1,X_1,...,X_N),$$
hence $H(F') = H(F)$. There exists $\bx = (x_0,...,x_N) \in \qbar^{N+1}$ so that $x_0 \neq 0$, and
$$F'(x_0,...,x_N) =  F(x_1/x_0,...,x_N/x_0) = 0.$$
Notice that
$$H(x_1/x_0,...,x_N/x_0) = H(x_1,...,x_N) \leq H(x_0,...,x_N) = H(\bx),$$
hence it is sufficient to prove that there exists a zero $\bz \in \qbar^{N+1}$ of $F'$ so that $z_0 \neq 0$ and $\bz$ is of bounded height. Notice that since the variable $X_0$ was introduced to homogenize $F$, we have $\deg(F) = \deg(F') = M$, and so $X_0 \nmid F'(X_0,...,X_N)$.

Write $\be'_0,...,\be'_N$ for the standard basis vectors in $\qbar^{N+1}$. First suppose that $F'(\be'_i) \neq 0$ for all $0 \leq i \leq N$, then by Theorem \ref{nonzero:coord} there exists $\bz \in \qbar^{N+1}$ satisfying (\ref{inh:bound}) with $\deg_K(\bz) \leq M$ such that $z_i \neq 0$ for each $0 \leq i \leq N$, and $F'(\bz)=0$, hence we are done. Next suppose that $F'(\be'_0)=F(\boldsymbol 0)=0$. Then let
$$G(X_1,...,X_N) = F'(X_1,X_1,...,X_N),$$
that is set $X_0=X_1$ in $F'$. Notice that for each $1 \leq i \leq N$, $G(\be_i) = F(\be_i) \neq 0$, and $H(G)=H(F')=H(F)$. Again, by Theorem \ref{nonzero:coord} there exists $\bz \in \qbar^N$ satisfying (\ref{inh:bound}) with $\deg_K(\bz) \leq M$ such that $z_i \neq 0$ for each $1 \leq i \leq N$, and $G(\bz) = F'(z_1,\bz) = 0$, and so we are done. Finally suppose that $F'(\be'_i)=0$ for some $1 \leq i \leq N$. Since  $X_0 \nmid F(X_0,...,X_N)$, we can write
$$F'(X_0, \dots, X_N) = G_1(X_1, \dots, X_N) + X_0 G_2(X_0, \dots, X_N),$$
where $G_1$ and $G_2$ are both non-zero homogeneous polynomials of degrees $M$ and $M-1$ respectively. Then $F'(\be'_i) = G_1(\be_i) = 0$, which means that the coefficient of the term $X_i^M$ in $G_1$ is zero, and hence it is zero in $F'$ and thus in $F$. This implies that $F(\be_i) = 0$ contradicting our original assumption. Hence $F'(\be'_i) \neq 0$ for every $1 \leq i \leq N$, and so we are done.
\endproof

In case $N=2$, the exponent in the bound of Corollary \ref{inhom} is best possible. Take
$$F(X_1,X_2) = X_1 - C X_2^M,$$
for some $0 \neq C \in K$. Then by the same argument as in the remark after the proof of Proposition \ref{bz:basic} every non-trivial zero of $F$ has height $ \geq O(H(F)^{1/M})$.
\bigskip

{\bf Acknowledgment.} I would like to thank Professors Paula Tretkoff and Jeff Vaaler for their helpful comments on the subject of this paper.

\bibliographystyle{plain}  
\bibliography{bezout}        

\end{document}